\documentclass{amsart}
\usepackage{amssymb,amsmath,latexsym}
\usepackage{amsthm}
\usepackage{fontenc}
\usepackage{amssymb}
\numberwithin{equation}{section}

\newtheorem{theorem}{Theorem}[section]

\begin{document}
\author{Alexander E Patkowski}
\title{On a symmetric $q$-series identity}

\maketitle
\begin{abstract}We prove an interesting symmetric $q$-series identity which generalizes a result due to Ramanujan. A proof that is analytic in nature is offered, and a bijective-type proof is also outlined.
\end{abstract}



\section{Introduction} In [1] Andrews gives a wonderful introduction of Ramanujan's ``Lost" notebook, and lists some interesting
identities contained therein. One of which is the following beautiful symmetric identity [1, eq.(1.5)], where if
\begin{equation}f(\alpha,\beta):=\frac{1}{1-\alpha}+\sum_{n\ge1}\frac{\beta^n}{(1-\alpha x^n)(1-\alpha x^{n-1}y)(1-\alpha x^{n-2}y^2)\cdots (1-\alpha y^n)},\end{equation}
then 
\begin{equation}f(\alpha,\beta)=f(\beta,\alpha).\end{equation}
The identity we present here is a refinement of the case where $x=q,$ and $y=q^2.$ Andrews provides an elegant bijective
proof of this identity in [1, pg.107] by taking the conjugate partition. We will also consider conjugate partitions in the third 
section, but will require a slightly different approach using a $2$-modular diagram (conceptually) to prove the following theorem bijectively. We first note
some notation which may be found in [6]. We put, throughout this paper, $(a)_n=(a;q)_{n}:=\prod_{0\le k\le n}(1-aq^{k}).$ Of course
the reader should note the infinite product that is obtained by passing the limit $n\rightarrow\infty,$ which we denote by $(a)_{\infty}.$
\begin{theorem} We have, for arbitrary $a,$ and $|b|<1,$ $|t|<1,$
\begin{equation}\sum_{n\ge0}\frac{(-abq^{n+1};q)_{n}t^n}{(bq^{n};q)_{n+1}}=\sum_{n\ge0}\frac{(-atq^{n+1};q)_{n}b^n}{(tq^{n};q)_{n+1}}.\end{equation}
\end{theorem}

\section{An analytic proof}
Since in [1] an analytic proof uses the binomial theorem, we decided to stay with our original use of a different $q$-polynomial identity. Namely,
we use the $q$-Pfaff-Saalsch$\ddot{u}$tz [6, pg.355, eq.(II.12)]
\begin{equation}\sum_{n\ge0}\frac{(a)_n(b)_n(q^{-N})_nq^n}{(c)_n(q)_n(q^{1-N}ab/c)_n}=\frac{(c/a)_n(c/b)_n}{(c)_n(c/ab)_N}.\end{equation}
The left side of (2.1) may be written
\begin{equation}\frac{(q)_N}{(c/(ab))_N} \sum_{n\ge0}\frac{(a)_n(b)_n(c/(ab))_{N-n}q^n}{(c)_n(q)_n(q)_{N-n}}(c/ab)^n.\end{equation}
Putting $c=bq$ in this identity we obtain 
\begin{equation} \sum_{n\ge0}\frac{(a)_n(q/a)_{N-n}}{(q)_n(q)_{N-n}(1-bq^n)}(q/a)^n=\frac{(bq/a)_N}{(b)_{N+1}}.\end{equation}
Now we may use (2.3) to compute the following:

$$\sum_{n\ge0}\frac{(abq^{n+1};q)_{n}t^n}{(bq^{n};q)_{n+1}}=\sum_{N\ge0}\sum_{n\ge0}\frac{(a)_n(q/a)_{N-n}(q/a)^nt^N}{(q)_n(1-bq^{N+n})(q)_{N-n}},$$

Shifting summation indices $N\rightarrow N+n$ gives,

$$\sum_{N\ge0}\sum_{n\ge0}\frac{(a)_n(q/a)_{N}(q/a)^nt^{N+n}}{(q)_n(1-bq^{N+2n})(q)_{N}}=\sum_{n\ge0}\frac{t^n(a)_n(q/a)^n}{(q)_n}\frac{(tq/a)_{\infty}}{(t)_{\infty}}\sum_{N\ge0}\frac{(t)_N}{(tq/a)_N}(bq^{2n})^N.$$

By a special case of [5, pg.18, eq.(16.3)] we compute this is equal to

$$\frac{(tq/a)_{\infty}}{(t)_{\infty} }\sum_{n\ge0}\frac{(t)_n}{(tq/a)_n}b^n \frac{(tq^{2n+1})_{\infty}}{(tq^{2n+1}/a)_{\infty}},$$

which may be simplified to the desired identity,

$$\sum_{n\ge0}\frac{(atq^{n+1};q)_{n}b^n}{(tq^{n};q)_{n+1}}.$$

\section{A Bijective Proof and some corollaries}
We first start with some standard notation on partitions, which can be found in [5, pg.37]. We write a partition $\pi$ to be a sequence which consists of nonnegative integers, say $(\pi_1, \pi_2, \cdots, \pi_m)$ where we say each $\pi_i$ for $1\le i \le m$ is a `part' with the largest $\pi_1,$ and smallest $\pi_m.$ The number of such parts is denoted $l(\pi),$ and the number of odd parts will be denoted $o(\pi).$ Since Guo obtained a similar symmetric $q$-series identity using partitions where odd parts do not repeat, we consider a similar approach. The main bijection appears to be due to R. Chapman in his proof of identities from [3] (see [4] amd [6] for more details). We will require an extra step in dealing with the inequality on parts that is in our identity, which is a key difference, however. We may replace $q$ by $q^2$ in (1.3) and then replace $a$ with $aq^{-1}$ to obtain that
\begin{equation}\sum_{n\ge0}\frac{(-abq^{2n+1};q^2)_{n}t^n}{(bq^{2n};q^2)_{n+1}}=\sum_{n\ge0}\frac{(-atq^{2n+1};q^2)_{n}b^n}{(tq^{2n};q^2)_{n+1}}.\end{equation}

Now, on the left hand side, $a$ keeps track of the number of odd parts, $b$ keeps track of the number of parts and $t$ keeps track of the largest part. It can then be seen that if we let $O$ be the set of partitions where odd parts do not repeat, we have that

\begin{equation}\sum_{\substack{\pi\in O\\ l(\pi)\le j\\ \pi_1\le 4M \\ \pi_m\ge 2M}}a^{o(\pi)}q^{|\pi|}=\sum_{\substack{\pi\in O\\ l(\pi)\le M\\ \pi_1\le 4j \\ \pi_m\ge 2j}}a^{o(\pi)}q^{|\pi|},\end{equation}
which has a similar resemblence (as is to be expected) to Guo's partition identity [7, Theorem 1.2]. The key difference is the inequality on the largest and smallest parts. This indeed causes a problem with using Chapman's bijection directly, but we have a simple solution to this. We consider the set of partitions $O^{*}$ where odd parts do not repeat, and the largest even appears each time a part appears. In our case we start with the left side of (3.2), and if parts are $\ge 2r,$ $r\in\mathbb{N},$ say, then $2r$ is removed from each part to appear as a separate part. This process ensures that in the new partition, $2r$ appears as a part each time a part appears. This set appears to us once we write out the factor in our $q$-series in the following way:
$$\frac{(1+abq^{2N+1})(1+abq^{2N+3})\cdots (1+abq^{2N-1+2N})}{(1-bq^{2N})(1-bq^{2+2N})\cdots (1-bq^{2N+2N})}.$$
If we consider pugging in natural numbers $N,$ we can conceptually see that we are rewriting each part into two different parts, one being $2N.$ So for example,
if $N=3:$ $7,9,11$ gets converted to $1+6,$ $3+6,$ and $5+6,$ ensuring each odd part is less than the largest even $(6, 6, 6, 5, 3, 1).$ Doing this to every part ensures that the number of parts in the new partition is less than or equal to double the number of parts in the previous partition. If we consider the set $O_{\pi_m\ge 2M, \pi_1\le 4M}$ to be the partitions from the left side of (3.2) then we have just created the mapping $\Gamma: O_{\pi_m\ge 2M, \pi_1\le 4M} \rightarrow O^{*}_{\pi_1\le 2M}.$ Once this process is completed we may take the conjugate partition of the $2$-modular diagram of this new partition to form a third partition, which satisfies the inequalities on the right side of (3.2). This second mapping is essentially the involution $\sigma$ of Chapman that was employed by Guo. We may apply $\sigma$ because $O^{*}$ is a partition where odd parts do not repeat. So since this preserves the sum of all the parts, and interchanges the number of parts and $\lceil \pi_1/2\rceil$ we see that the composition $\sigma \circ \Gamma: O_{\pi_m\ge 2M, \pi_1\le 4M}\rightarrow O_{\pi_m\ge 2j, \pi_1\le 4j}$ gives the theorem since this is 1-to-1 (both $\sigma$ and $\Gamma$ are reversible). The inequality in the codomain is obtained because the smallest part in the codomain is twice the number of largest parts in the domain (the number of largest parts from $O^{*}).$ In summation notation this may be written:
\begin{equation}\sum_{\substack{\pi\in O\\ l(\pi)\le j\\ \pi_1\le 4M \\ \pi_m\ge 2M}}a^{o(\pi)}q^{|\pi|}=\sum_{\substack{\lambda\in O^{*}\\ l(\lambda)\le 2j\\ \lambda_1\le 2M}}a^{o(\lambda)}q^{|\lambda|}=\sum_{\substack{\pi\in O\\ l(\pi)\le M\\ \pi_1\le 4j \\ \pi_m\ge 2j}}a^{o(\pi)}q^{|\pi|}.\end{equation}

Next we offer two examples. We denote $\pi^{*}$ to be a partition taken from the set $O^{*}$ and $\pi'$ to be the result of taking the conjugate of a $2$-modular diagram.
\\*

{\bf Example 1.} \it We consider a partition from the set $O$ where parts are $\ge10$ and $\le20.$ Take $\pi=(20,17,13),$ then $\pi^{*}=(10,10,10,10, 7,3),$ and so $\pi'=(12, 11, 9, 8).$ Notice $l(\pi)=3,$ and so in $\pi'$ we have that $\pi'_{m}\ge6,$ $\pi'_{1}\le12.$ (We may have done this a little bit differently if we had picked $12$ to be the largest part instead of $10,$ which would have resulted in $\pi^{*}=(12,12,12,8,5,1),$ and $\pi'=(11,10,9,8,6,6).$)
\\*
\rm

{\bf Example 2.} \it We consider a partition from the set $O$ where parts are $\ge10$ and $\le20.$ Take $\pi=(20,13,12,12,10),$ then $\pi^{*}=(10,10,10,10,10,10,3,2,2),$ and so $\pi'=(18,13,12,12,12).$ Notice $l(\pi)=5,$ and so in $\pi'$ we have that $\pi'_{m}\ge10,$ $\pi'_{1}\le20.$
\\*
\rm

In closing we present a few charming identities that may be obtained by employing the work from [3].

\begin{theorem} We have,
\begin{equation} \sum_{N\ge0}\left(1-\frac{(abq^{N+1})_N}{(bq^N)_N}\right)=-\sum_{n\ge1}\frac{(aq^{n+1})_n}{(q^n)_{n+1}}b^n,\end{equation}
\begin{equation} \sum_{N\ge0}\left(1-(bq^{N+1})_N\right)=-\sum_{n\ge1}\frac{(-b)^nq^{n(3n+1)/2}}{(q^n)_{n+1}}.\end{equation}
\end{theorem}

1390 Bumps River Rd. \\*
Centerville, MA
02632 \\*
USA \\*
E-mail: alexpatk@hotmail.com
\end{document}